\documentclass[final]{siamltex}

\usepackage{amsmath,amsfonts,graphicx}

 \newtheorem{thm}{Theorem}[section]
 
 \newtheorem{lem}[thm]{Lemma}
 \newtheorem{prop}[thm]{Proposition}
 \newtheorem{assumption}[thm]{Assumption}
 \newtheorem{rem}[thm]{Remark}
 
 \numberwithin{equation}{section}
\def\F{\mathcal{F}}
\def\R{\mathbb{R}}
\def\N{\mathbb{N}}
\def\B{\mathbf{B}}
\def\K{\mathbf{K}}
\def\Q{\mathrm{Q}}
\def\C{\mathbf{C}}
\def\P{\mathrm{P}}

\def\L{\mathrm{L}}

\def\f{\mathbf{f}}

\def\y{\mathbf{y}}

\def\b{\mathcal{B}}

\def\vol{{\rm vol}\,}

\def\figurewidth{0.54}

\title{Approximate volume and integration for basic semi-algebraic sets\thanks{This work was completed with the support of 
the (French) ANR grant NT05-3-41612. The first author was also supported by
the research program No. MSM6840770038
of the Czech Ministry of Education and Grant  102/08/0186 of
the Grant Agency of the Czech Republic. }}

\author{D. Henrion\thanks{LAAS-CNRS, University of Toulouse,
7 avenue du Colonel Roche, 31077 Toulouse, France, and
Faculty of Electrical Engineering, Czech Technical University in
Prague, Technick\'a 2, 16627 Prague, Czech Republic ({\tt henrion@laas.fr}).}
\and J. B. Lasserre\thanks{
LAAS-CNRS, University of Toulouse,
7 avenue du Colonel Roche, 31077 Toulouse, France, and  Institute of Mathematics,
University of Toulouse ({\tt lasserre@laas.fr}).}
\and C. Savorgnan\thanks{
Formerly with LAAS-CNRS, University of Toulouse,
7 avenue du Colonel Roche,
31077 Toulouse,
France. Now with the Department of Electrical Engineering,
Katholieke Universiteit Leuven, Belgium ({\tt carlo.savorgnan@gmail.com}).}}

\begin{document}

\maketitle

\begin{abstract}
Given a basic compact semi-algebraic set $\K\subset\R^n$,
we introduce a methodology that generates
a sequence converging to the volume of $\K$.
This sequence is obtained from optimal values of a hierarchy of either semidefinite 
or linear  programs.
Not only the volume but also every finite vector of moments of the probability measure 
that is uniformly distributed on $\K$ can be approximated as closely as desired,
and so permits to approximate the integral on $\K$ of any given polynomial;
extension to integration against some weight functions is also provided.
Finally, some numerical issues associated with the algorithms involved are 
briefly discussed.
\end{abstract}

\begin{keywords}Computational geometry; volume; integration;
$\K$-moment problem; semidefinite programming
\end{keywords}

\begin{AMS}
14P10, 11E25, 12D15, 90C25
\end{AMS}

\pagestyle{myheadings}
\thispagestyle{plain}
\markboth{D. Henrion, J.B. Lasserre and C. Savorgnan}{Approximate volume and integration}

\section{Introduction}~
Computing the volume and/or integrating on a subset $\K\subset\R^n$
is a challenging problem with many important applications. 
One possibility is to use basic Monte Carlo techniques that generate points
uniformly in a box containing $\K$ and then count the proportion of 
points falling into $\K$. To the best of our knowledge, all other approximate (deterministic or randomized) or exact techniques deal with polytopes or convex bodies only. Similarly,
powerful cubature formulas exist for numerical integration against a weight function 
on simple sets (like e.g. simplex, box),  but not for arbitrary semi-algebraic sets.

The purpose of this paper is to introduce a deterministic technique
that potentially applies to any basic compact semi-algebraic set $\K\subset\R^n$.
It is deterministic (no randomization) and differs from previous ones in the literature
essentially dedicated to {\it convex} bodies (and more particularly, convex polytopes). Indeed, one treats the original problem as an infinite dimensional optimization (and even linear programming (LP)) problem
whose unknown is the Lebesgue measure
on $\K$. 
Next, by restricting to finitely many of its moments, and using a 
certain characterization on the $\K$-moment problem, one ends up 
in solving a hierarchy of semidefinite programming (SDP) problems whose size 
is parametrized by the number of moments considered;
the dual LP has a simple interpretation 
and from this viewpoint, convexity of $\K$ does not help much.
For a certain choice 
of the criterion to optimize, one obtains a monotone non increasing sequence of upper bounds on 
the volume of $\K$. Convergence to the exact value
invokes results on the $\K$-moment problem by Putinar \cite{putinar}.
Importantly, there is {\it no} convexity and not even connectedness assumption on $\K$,
as this plays no role in the $\K$-moment problem.
Alternatively, using a different 
characterization of the $\K$-moment problem due to Krivine \cite{krivine}, one 
may solve a hierarchy of LP (instead of SDP) problems
whose size is also parametrized by the number of moments.
Our contribution is a new addition to the already very long list
of applications of the moment approach (some of them described in
e.g. Landau \cite{landau} and Lasserre \cite{lasserre3}) and
semidefinite programming \cite{boyd}.  In principle,
the method also permits to approximate any finite number of moments
of the uniform distribution on $\K$, and so provides a means to approximate the integral
of a polynomial on $\K$. Extension to integration against a weight function is also proposed.\\

\noindent
{\bf Background.} Computing or even approximating 
the volume of a convex body is hard theoretically and in practice as well.
Even if $\K\subset\R^n$ is a convex polytope, exact computation of its volume or
integration over $\K$ is a computational challenge. Computational complexity of these problems
is discussed in e.g. Bollob\'as \cite{bollobas} and Dyer and Frieze \cite{dyer1}.
Any deterministic algorithm with polynomial time
complexity that would compute an upper bound $\overline{\vol(\K)}$ and a lower bound
$\underline{\vol(\K)}$ on $\vol(\K)$
cannot yield an upper bound on the ratio $\overline{\vol(\K)}/\underline{\vol(\K)}$ better
than polynomial in the dimension $n$. Methods for exact volume computation 
use either triangulations or simplicial decompositions
depending on whether the polytope has a half-space description or
a vertex description. See e.g. Cohen and Hickey, \cite{cohen}, Lasserre \cite{lasserrejota},
Lawrence \cite{lawrence} and see B\"ueler et al. \cite{bueler} for a comparison.
Another set of methods which use generating functions are described in e.g. Barvinok \cite{barvinok}
and Lasserre and Zeron \cite{laszeron}. Concerning integration on simple sets 
(e.g. simplex, box) via cubature formulas, the interested reader is referred to
Gautschi \cite{gaut1,gaut2} and Trefethren \cite{tref}.

A {\it convex body} $\K\subset\R^n$ is a compact convex subset with nonempty interior. A strong separation {\it oracle} answers either $x\in \K$ or $x\not\in \K$, and in the latter case produces a hyperplane separating $x$ from $\K$. A negative result states that for every polynomial-time 
algorithm for computing the volume of a convex body $\K\subset\R^n$
given by a well-guaranteed separation oracle, there is a constant $c>0$ such that
$\overline{\vol(\K)}/\underline{\vol(\K)}\leq(cn/\log n)^n$
cannot be guaranteed for $n\geq2$.
However, Lov\'asz \cite{lovasz} proved that there is a polynomial-time algorithm that produces
$\overline{\vol(\K)}$ and $\underline{\vol(\K)}$ satisfying
$\overline{\vol(\K)}/\underline{\vol(\K)}\leq n^n\,(n+1)^{n/2}$, whereas
Elekes \cite{elekes} proved that for $0<\epsilon<2$ there is no polynomial-time algorithm
that produces $\overline{\vol(\K)}$ and $\underline{\vol(\K)}$ satisfying
$\overline{\vol(\K)}/\underline{\vol(\K)}\leq (2-\epsilon)^n$.

If one accepts randomized algorithms that fail with small probability, 
then the situation is more favorable.
Indeed, the celebrated Dyer, Frieze and Kanan probabilistic 
approximation algorithm \cite{dyer2}
computes the volume to fixed arbitrary relative precision $\epsilon$,
in time polynomial in $\epsilon^{-1}$. The latter algorithm
uses approximation schemes based on rapidly mixing Markov chains and 
isoperimetric inequalities.  See also  hit-and-run algorithms for sampling points
according to a given distribution, described in e.g. 
Belisle \cite{belisle}, Belisle et al. \cite{smith2}, and Smith \cite{smith1}\\

\noindent
{\bf Contribution.}
This paper is concerned with computing (or rather approximating) the volume of a
compact basic semi-algebraic set $\K\subset\R^n$ defined by
\begin{equation}
\label{setk}
\K\,:=\,\{\:x\in\R^n\::\: g_j(x)\,\geq\,0,\quad j=1,\ldots,m\:\}
\end{equation}
for some polynomials $(g_j)_{j=1}^m\subset\R[x]$. Hence $\K$ is possibly 
non-convex and non-connected. Therefore, in view of
the above discussion, this is quite a challenging problem.

(a) We present a numerical scheme that depends on a parameter $p$,
a polynomial that is nonnegative on $\K$ (e.g. $p\equiv 1$).
For each parameter $p$, it provides converging approximations of moments
of the measure uniformly supported on $\K$ (with mass
equal to $\vol(\K)$). For the choice $p\equiv 1$ one obtains
a {\it monotone} non-increasing sequence of upper bounds that converges to
$\vol(\K)$. 

(b) The way we see the problem dates back to the 19th century pioneer work in the 
one-dimensional case by
Chebyshev \cite{cebicev}, Markov \cite{markov} and Stieltjes \cite{stieltjes}, 
where given $n$ moments $s_k=\int_a^b t^kf(t)dt$, $k=0,\ldots,n-1$,
and $a <c<d<b$, one wishes to approximate the integral
$\int_c^d f(t)dt$ and analyzes asymptotics as $n\to\infty$;
characterizing feasible sequence $(s_k)$ is referred to as the 
Markov moment problem (and $L$-moment problem if in addition
one requires $0\leq f\leq L$ for some scalar $L$).
For an historical account on this problem as well as other developments,
the interested reader is referred to e.g. Krein \cite{krein},
Krein and Nuldelman \cite{krein2}, Karlin and Studden \cite{karlin} and Putinar \cite{putinar2}.

Our method combines a simple idea, easy to describe, with
relatively recent powerful results on the $\K$-moment problem described in e.g. 
\cite{claus,schmudgen,markus}.
It only requires knowledge of a set $\B$ (containing $\K$) 
simple enough so that the moments of the Lebesgue measure on $\B$
can be obtained easily. For instance $\B := \{x \in {\mathbb R}^n \: :\:
\|x\|_p\leq a\}$ with $p=2$ (the scaled $n$-dimensional ball) or
$p=\infty$ (the scaled $n$-dimensional box) and $a \in \mathbb R$
a given constant.
Then computing $\vol(\K)$ is equivalent to computing
the mass of the Borel measure $\mu$ which is
the restriction to $\K$ of the Lebesgue measure
on $\B$.  This in turn is translated
into an infinite dimensional LP problem $\P$ 
with parameter $p$ (some polynomial nonnegative on $\K$) and 
with the Borel measure $\mu$ as unknown.
Then, from certain results on the $\K$-moment problem
and its dual theory of the representation
of polynomials positive on $\K$, problem $\P$ can be approximated by an appropriate
hierarchy of semidefinite programs (SDP) whose size depends on
the number $d$ of moments of $\mu$ considered.
One obtains approximations of the moments of $\mu$ which converge to the exact value as $d\to\infty$.
For the choice $p\equiv 1$ of the parameter $p$, one even 
obtains an non-increasing sequence of upper bounds 
converging to $\vol(\K)$.
Asymptotic convergence is ensured by invoking results of Putinar
\cite{putinar} on the $\K$-moment problem.
Alternatively, one may replace the SDP hierarchy with an LP hierarchy
and now invoke results of Krivine \cite{krivine} for convergence. 

Interestingly, the dual of each SDP relaxation defines
a strenghtening of $\P^*$, the LP dual of $\P$,
and highlights why the problem of computing the volume is difficult.
Indeed, one has to approximate from above the function
$f$ ($=p$ on $\K$ and $0$ on $\B\setminus\K$)
by a polynomial $h$ of bounded degree,
so as to minimize the integral 
$\int_{\B}(h-f)dx$. From this viewpoint, convexity of $\K$ plays no particular 
role and so, does not help much.

(c) Let $d\in\N$ be fixed, arbitrary. One 
obtains an approximation of the moments of degree up to $d$ of the measure 
$\mu$ on $\K$,
as closely as desired.
Therefore, this technique also provides a sequence of approximations 
that converges to $\int_\K q dx$ for  {\it any} polynomial $q$
of degree at most $d$ (in contrast, Monte Carlo simulation
is for a given $q$). Finally, we also propose a similar approximation scheme for
integrating a polynomial on $\K$ against a nonnegative weight function $w(x)$. The only required data 
are moments of the measure $d\nu=wdx$ on a simple set $\B$ (e.g. box or simplex) containing $\K$,
which can be obtained by usual cubature formulas for integration.\\

On the practical side, at each step $d$ of the hierarchy, 
the computational workload is that of solving an SDP
problem of
increasing size. In principle, this can be done in time polynomial in the input size of the SDP problem, at
given relative accuracy.
However, in view of the present status and limitations of current SDP solvers, 
so far the method is restricted to problems of small dimension $n$
if one wishes to obtain good approximations. The alternative 
LP hierarchy might be preferable for larger size problems,
even if proved to be less efficient when used in other contexts where the moment approach applies,
see e.g. \cite{lasserre2,lasserre-prieto}. 

Preliminary results on simple problems for which $\vol(\K)$ is known show
that indeed convexity plays no particular role. In addition, 
as for interpolation problems, the choice of
the basis of polynomials is crucial from the viewpoint of numerical precision. This is 
illustrated on a trival example on the real line where, as expected,  the basis of Chebyshev
polynomials is far better than the usual monomial basis. In fact,
it is conjectured that trigonometric polynomials would be probably the best choice.
Finally, the choice of the parameter $p$ is also very important and unfortunately, 
the choice of $p\equiv 1$ which guarantees a monotone convergence to $\vol(\K)$ is 
not the best choice at all.  Best results are obtained when $p$ is negative outside $\K$.

So far, for convex polytopes, this method is certainly not competitive 
with exact specific methods as those described
in e.g. \cite{bueler}.
It rather should be viewed as a relatively simple deterministic methodology 
that applies to a very general context for which even getting good bounds on 
$\vol(\K)$ is very difficult, and for which the only 
alternative presently available seems to be brute force Monte Carlo.

\section{Notation, definitions and preliminary results}

Let $\R[x]$ be the ring of real polynomials in the variables $x=(x_1,\ldots,x_n)$,
and let $\Sigma^2[x]\subset\R[x]$ be the subset of sums of squares (SOS) polynomials. 
Denote $\R[x]_d\subset\R[x]$ be the set of
polynomials of degree at most $d$, which forms a vector space of dimension
$s(d)={n+d\choose d}$.
If $f\in\R[x]_d$, write
$f(x)=\sum_{\alpha\in\N^n}f_\alpha x^\alpha$ in
the usual canonical basis $(x^\alpha)$, and
denote by $\f=(f_\alpha)\in\R^{s(d)}$ its vector of coefficients.  
Similarly, denote by $\Sigma^2[x]_d\subset\Sigma^2[x]$ the subset of
SOS polynomials of degree at most $2d$.
\vspace{0.2cm}

\noindent
{\bf Moment matrix.} Let $\y=(y_\alpha)$ be a sequence indexed in the canonical basis
$(x^\alpha)$ of $\R[x]$, let $L_\y:\R[x]\to\R$ be the linear functional
\[f\quad (=\sum_{\alpha}f_{\alpha}\,x^\alpha)\quad\mapsto\quad
L_\y(f)\,=\,\sum_{\alpha}f_{\alpha}\,y_{\alpha},\]
and let $M_d(\y)$ be the symmetric matrix with rows and columns indexed in 
the canonical basis $(x^\alpha)$, and defined by:
\[M_d(\y)(\alpha,\beta)\,:=\,L_\y(x^{\alpha+\beta})\,=\,y_{\alpha+\beta},\]
for every $\alpha,\beta\in\N^n_d:=\{\alpha\in\N^n\::\:\vert \alpha\vert \:(=\sum_i\alpha_i)\leq d\}$.

A sequence $\y=(y_\alpha)$ is said to have a {\it representing} finite Borel measure $\mu$ if
$y_\alpha=\int x^\alpha d\mu$ for every $\alpha\in\N^n$. A necessary 
(but not sufficient) condition is that $M_d(\y)\succeq0$ for every $d\in\N$. However,
 if in addition, $\vert y_\alpha\vert\leq M$ for some M and for every $\alpha\in\N^n$, then $\y$
  has a representing measure on $[-1,1]^n$.

\vspace{0.2cm}

\noindent
{\bf Localizing matrix.} Similarly, with $\y=(y_{\alpha})$
and $g\in\R[x]$ written as
\[x\mapsto g(x)\,=\,\sum_{\gamma\in\N^n}g_{\gamma}\,x^\gamma,\]
let $M_d(g\,\y)$ be the symmetric matrix with rows and columns indexed in 
the canonical basis $(x^\alpha)$, and defined by:
\[M_d(g\,\y)(\alpha,\beta)\,:=\,L_\y\left(g(x)\,x^{\alpha+\beta}\right)\,=\,\sum_{\gamma}g_{\gamma}\,
y_{\alpha+\beta+\gamma},\]
for every $\alpha,\beta\in\N^n_d$. A necessary (but not sufficient) condition for $\y$ to have a representing measure 
with support contained in the level set $\{x\,:\, g(x)\geq0\}$ is that $M_d(g\,\y)\succeq0$ for every $d\in\N$.
\vspace{0.2cm}

\subsection{Moment conditions and representation theorems}

The following results from the $\K$-moment problem and its dual theory of
polynomials positive on $\K$ provide the rationale behind
the hierarchy of SDP relaxations introduced in \cite{lasserre1}, and 
potential applications in many different contexts. See e.g. \cite{lasserre3} and the many references therein. \\

\noindent
{\bf SOS-based representations.}
Let $Q(g)\subset\R[x]$ be the quadratic module generated by polynomials $(g_j)_{j=1}^m\subset\R[x]$, that is,
\begin{equation}
\label{qg}
Q(g)\,:=\,\left\{\sigma_0+\sum_{j=1}^m\sigma_j\,g_j\::
\quad(\sigma_j)_{j=1}^m\subset\Sigma^2[x]\:\right\}.\end{equation}
\begin{assumption}
\label{assput}
The set $\K\subset\R^n$  in (\ref{setk}) is compact and the quadratic polynomial $x\mapsto a^2-\Vert x\Vert^2$ belongs to $Q(g)$ for some given constant $a \in {\mathbb R}$.
\end{assumption}
\begin{thm}[Putinar's Positivstellensatz \cite{putinar}]
\label{thput}
Let Assumption \ref{assput} hold.

{\rm (a)} If $f\in\R[x]$ is strictly positive on $\K$, then $f\in Q(g)$. That is:
\begin{equation}
\label{putinarrep}
f\,=\,\sigma_0+\sum_{j=1}^m\sigma_j\, g_j,\end{equation}
for some SOS polynomials $(\sigma_j)_{j=1}^m\subset\Sigma^2[x]$. 

{\rm (b)} If  $\y=(y_\alpha)$ is such that for every $d\in\N$, 
\begin{equation}
\label{putinarrep2}
M_d(\y)\succeq0;\quad M_d(g_j\y)\succeq0,\qquad j=1,\ldots,m, 
\end{equation}
then $\y$ has a representing finite Borel measure $\mu$ supported on $\K$.
\end{thm}

Given $f\in\R[x]$, or $\y=(y_\alpha)\subset\R$, checking 
whether (\ref{putinarrep}) holds for  SOS $(\sigma_j)\subset\Sigma^2[x]$ with a priori
bounded degree, or checking whether (\ref{putinarrep2}) holds with $d$ fixed, 
reduces to solving an SDP.\\

\noindent
{\bf Another type of representation.}
Let $\K\subseteq\B$ be as in (\ref{setk}) and assume for simplicity that the $g_j$s have been scaled to satisfy
$0\leq g_j\leq 1$ on $\K$, for every $j=1,\ldots,m$.  In addition, assume that the family of polynomials
$(1,g_1,\ldots,g_m)$ generates the algebra $\R[x]$.   For every
$\alpha\in\N^m$, let $g^\alpha$ and $(1-g)^\beta$ denote the polynomials
\[x\mapsto \quad g(x)^\alpha :=g_1(x)^{\alpha_1}\cdots g_m(x)^{\alpha_m},\]
and 
\[x\mapsto\quad(1-g(x))^\beta :=(1-g_1(x))^{\beta_1}\cdots (1-g_m(x))^{\beta_m}.\]
the following result is due to Krivine \cite{krivine} but is explicit in e.g. Vasilescu \cite{vasi}.

\begin{thm}
\label{thkrivine}

{\rm (a)} If $f\in\R[x]$ is strictly positive on $\K$, then
\begin{equation}
\label{krivinerep}
f\,=\,\sum_{\alpha,\beta\in\N^m}\,c_{\alpha\beta}\,g^\alpha\,(1-g)^\beta
\end{equation}
for finitely many nonnegative scalars $(c_{\alpha\beta})\subset\R_+$.

{\rm (b)} If  $\y=(y_\alpha)$ is such that
\begin{equation}
\label{krivinerep2}
L_\y(g^\alpha\,(1-g)^\beta)\geq\,0,
\end{equation}
for every $\alpha,\beta\in\N^m$,
then $\y$ has a representing finite Borel measure $\mu$ supported on $\K$.
\end{thm}
\vspace{0.2cm}

Theorem \ref{thkrivine} extends the well-known Hausdorff moment conditions on
the hyper cube $[0,1]^n$, as well as Handelman representation \cite{handelman} for convex polytopes
$\K\subset\R^n$.
Observe that checking whether (\ref{krivinerep}), resp. (\ref{krivinerep2}), holds
with $\alpha,\beta$ bounded a priori, reduces to solving an LP
in the variables $(c_{\alpha\beta})$, resp. $(y_\alpha)$.

\subsection{A preliminary result}

Given any two measures $\mu_1,\mu_2$ on a Borel $\sigma$-algebra
$\b$, the notation $\mu_1\leq\mu_2$ means $\mu_1(\C)\leq\mu_2(\C)$ for every $\C\in\b$. 

\begin{lem}
\label{ll}
Let Assumption \ref{assput} hold and let $\y_1=({y_1}_\alpha)$ and $\y_2=({y_2}_\alpha)$
be two moment sequences with respective
representing measures $\mu_1$ and $\mu_2$ on $\K$.
If \[M_d(\y_2-\y_1)\succeq0\,;\qquad M_d(g_j\,(\y_2-\y_1))\succeq 0,\quad j=1,\ldots,m,\]
for every $d\in\N$, then $\mu_1\leq\mu_2$.
\end{lem}
\begin{proof}
As $M_d(\y_2-\y_1)\succeq0$ and $M_d(g_j\,(\y_2-\y_1))\succeq 0$
for $j=1,\ldots,m$ and $d\in\N$, by Theorem \ref{thput}, the sequence
$\y_0:=\y_2-\y_1$ has a representing Borel measure $\mu_0$ on $\K$. From
${y_0}_\alpha+{y_1}_\alpha={y_2}_\alpha$ for every $\alpha\in\N^n$,
we conclude that \[\int x^\alpha\,d\mu_0+\int x^\alpha\,d\mu_1\,=\,\int x^\alpha\,d\mu_2,
\qquad\forall \alpha\in\N^n,\]
and as $\K$ is compact, by the Stone-Weierstrass theorem,
\[\int f\,d\mu_0+\int f\,d\mu_1\,=\,\int f\,d\mu_2 \]
for every continuous function $f$ on $\K$, which in turn implies
$\mu_0+\mu_1=\mu_2$, i.e., the desired result $\mu_1\leq\mu_2$.
\end{proof}

\section{Main result}

We first introduce an infinite-dimensional LP problem
$\P$ whose unique optimal solution is the restriction $\mu$ of 
the normalized Lebesgue measure on $\B$ (hence with 
$\mu(\K)=\vol(\K)/2^n$) and whose dual
has a clear interpretation. We then define a hierarchy of SDP problems
(alternatively, a hierarchy of LP problems)
to approximate any finite sequence of moments of $\mu$, as closely as desired.

\subsection{An infinite-dimensional linear program $\P$}
After possibly some normalization of the defining polynomials,
assume with no loss of generality that
$\K\subset\B\subseteq [-1,1]^n$ with $\B$ a set over which integration w.r.t. the Lebesgue measure
is easy. For instance, $\B$ is the box $[-1,1]^n$ or
$\B$ is the euclidean unit ball. 

Let $\b$ be the Borel $\sigma$-algebra of Borel subsets of $\B$, and 
let $\mu_2$ be the Lebesgue measure on $\B$, normalized so that
$2^n\mu_2(\B)={\rm vol}(\B)$. Therefore, if $\vol(\C)$ denotes
the $n$-dimensional volume of $\C\in\b$, then
$\mu_2(\C)=\vol(\C)/2^n$ for every $\C\in\b$.

Also, the notation $\mu_1\ll\mu_2$ means that $\mu_1$ is absolutely continuous w.r.t.
$\mu_2$, and $L_1(\mu_2)$ is the set of all functions integrable w.r.t. $\mu_2$.
By the Radon-Nikodym theorem, there exists a nonnegative
measurable function $f\in L_1(\mu_2)$ such that $\mu_1(\C)=\int_\C fd\mu_2$ for every $\C\in\b$,
and $f$ is called the Radon-Nikodym derivative of $\mu_1$ w.r.t. $\mu_2$.
In particular, $\mu_1\leq\mu_2$ obviously implies $\mu_1\ll\mu_2$.
For $\K\in\b$, let $M(\K)$ be the set of finite Borel measures on $\K$.

\begin{thm}
\label{th-mesures}
Let $\K\in\mathcal{B}$ with $\K\subseteq\B$ and let $p\in\R[x]$ be positive
almost everywhere on $\K$. Consider the following infinite-dimensional LP problem:
\begin{equation}
\label{defp}
\P:\quad \displaystyle\sup_{\mu_1}\:
\{\,\int p\,d\mu_1\::\:\mu_1\leq\mu_2;\quad\mu_1\in M(\K)\:\}
\end{equation}
with optimal value denoted $\sup\P$ (and $\max\P$ if the supremum is achieved).

Then the restriction $\mu_1^*$ of $\mu_2$ to $\K$
is the unique optimal solution of $\P$ and
$\max\P=\int pd\mu_1^*=\int_\K pd\mu_2$.
In particular, if $p\equiv 1$ then $\max\P=\vol(\K)/2^n$.
\end{thm}
\begin{proof}
Let $\mu_1^*$ be the restriction of $\mu_2$ to $\K$
(i.e. $\mu_1^*(\C)=\mu_2(\C\cap \K)$, $\forall\,\C\in\b$).
Observe that $\mu_1^*$ is a feasible solution of $\P$.
Next, let $\mu_1$ be any feasible solution of $\P$. 
As $\mu_1\leq\mu_2$ then 
\[\mu_1(\C\cap \K)\,\leq\mu_2(\C\cap \K)\,=\,\mu_1^*(\C\cap\K),\qquad \forall\,\C\in\b,\]
and so, $\mu_1\leq\mu_1^*$ because $\mu_1$ and $\mu_1^*$ are supported on $\K$.
Therefore, as $p\geq0$ on $\K$,
$\int pd\mu_1\leq\int pd\mu_1^*$ which proves that $\mu_1^*$ is an optimal solution of $\P$.

Next suppose that $\mu_1\neq\mu_1^*$ is another optimal solution of $\P$.
As $\mu_1\leq\mu_1^*$ then $\mu_1\ll\mu_1^*$ and so, by the Radon-Nikodym theorem,
there exists a nonnegative measurable function
$f\in L_1(\mu_1^*)$ such that 
\[\mu_1(\C)\,=\,\int_{\C}d\mu_1\,=\,\int_{\C} f(x)\,d\mu_1^*(x),\qquad \forall\,\C\in\b\cap\K.\]
Next, as $\mu_1\leq\mu_1^*$, $\mu_1^*-\mu_1=:\mu_0$ is a finite Borel measure on $\K$
which satisfies
\[0\,\leq\,\mu_0(\C)\,=\,\int_{\C}(1-f(x))\,d\mu_1^*(x),\qquad \forall\, {\C}\in\b\cap\K,\]
and so $1\geq f(x)$ for almost all $x\in\K$. 
But then, since $\int pd\mu_1=\int pd\mu_1^*$,
\[0\,=\,\int pd\mu_0\,=\,\int_\K p(x)(1-f(x))\,d\mu_1^*(x),\]
which (recalling $p>0$ almost everywhere on $\K$)
implies that $f(x)=1$ for almost-all $x\in\K$. And so
$\mu_1=\mu_1^*$.
\end{proof}

\subsection{The dual of $\P$}
Let $\F$ be the Banach space of continuous functions on $\B$ (equipped with the sup norm) and $\F_+$ its positive cone, i.e., 
the elements $f\in\F$ which are nonnegative on $\B$. The dual of $\P$ reads:
\begin{equation}
\label{lpdual}
\P^*:\quad\displaystyle\inf_{f\in\F_+}\:\{\int f\,d\mu_2\::\:f\geq p\mbox{ on }\K\}
\end{equation}
with optimal value denoted $\inf\P^*$ ($\min\P^*$ is the infimum is achieved).

Hence, a minimizing sequence of $\P^*$ aims at approximating from above
the function $f$ ($=p$ on $\K$ and $0$ on $\B\setminus\K$)
by a sequence $(f_\ell)$ of continuous functions
so as to minimize $\int f_\ell d\mu_2$.

Let $x\mapsto d(x,\K)$ be the euclidean distance to the set $\K$ and with $\epsilon_\ell>0$,
let $\K_\ell:=\{x\in\B\::\:d(x,\K)< \epsilon_\ell\}$ be an open bounded outer approximation of $\K$,
so that $\B\setminus\K_\ell$ is closed (hence compact) with $\epsilon_\ell\to 0$ as $\ell\to\infty$. 
By Urysohn's Lemma \cite[A4.2, p. 379]{ash}, there exists a sequence 
$(f_\ell)\subset\F_+$ such that
$0\leq f_\ell\leq 1$ on $\B$, $f_\ell=0$ on $\B\setminus \K_\ell$,
and $f_\ell=1$ on $\K$. Therefore,
\[\int f_\ell\,d\mu_2\,=\, \vol(\K)/2^n+\int_{\K_\ell\setminus \K}f_\ell d\mu_2,\]
and so $\int f_\ell d\mu_2\to\vol(\K)/2^n$ as $\ell\to\infty$. Hence, for the choice of the parameter $p\equiv 1$,
$\vol(\K)/2^n$ is the optimal value of both $\P$ and $\P^*$. 

\subsection{A hierarchy of semidefinite relaxations for computing the volume of $\K$}
Let $\y_2=({y_2}_\alpha)$ be the sequence of all moments of $\mu_2$. For example,
if $\B=[-1,1]^n$, then
\[{y_2}_\alpha = 2^{-n}\prod_{j=1}^n \left(\frac{2((1+\alpha_j)\bmod 2)}{1+\alpha_j}
\right)
,\qquad \forall\,\alpha\in\N^n.\]
Let $\K$ be a compact semi-algebraic set
as in (\ref{setk}) and let $r_j=\lceil ({\rm deg}\,g_j)/2\rceil$, $j=1,\ldots,m$.
Let $p\in\R[x]$ be a given polynomial 
positive almost everywhere on $\K$, and let $r_0:=\lceil ({\rm deg}\,p)/2\rceil$. For $d\geq \max_jr_j$, consider the following semidefinite program:
\begin{equation}
\label{relaxdef}
\Q_d:\quad\left\{\begin{array}{lll}
\displaystyle\sup_{\y_1} &L_{\y_1}(p)&\\
\mbox{s.t.}&M_d(\y_1)&\succeq\,0\\
&M_d(\y_2-\y_1)&\succeq\,0\\
&M_{d-r_j}(g_j\,\y_1)&\succeq\,0,\qquad j=1,\ldots,m
\end{array}\right.
\end{equation}
with optimal value denoted $\sup\Q_d$ (and $\max\Q_d$ if the supremum is achieved).

Observe that $\sup\Q_d\geq\max\P$ for every $d$. Indeed, the sequence $\y_1^*$ of moments of 
the Borel measure $\mu_1^*$ (restriction of $\mu_2$ to $\K$ and unique optimal solution of $\P$)
is a feasible solution of $\Q_d$ for every $d$.

\begin{thm}\label{th-sdp}
Let Assumption \ref{assput} hold and consider the hierarchy of
semidefinite programs $(\Q_d)$ in (\ref{relaxdef}). Then:

{\rm (a)} $\Q_d$ has an optimal solution (i.e. $\sup\Q_d=\max\Q_d$)
and 
\[\max\Q_d\downarrow \int_\K p\,d\mu_2,\qquad\mbox{as $d\to\infty$.}\]

{\rm (b)} Let ${\y_1}^d=({y_1}^d_\alpha)$ be an optimal solution of $\Q_d$, then
\begin{equation}
\label{conv}
\lim_{d\to\infty}
{y_1}^d_\alpha\,=\, \int_\K x^\alpha\,d\mu_2,\qquad \forall \alpha\in\N^n.
\end{equation}
\end{thm}
\begin{proof}
(a) and (b). Recall that $\B\subseteq [-1,1]^n$. 
By definition of $\mu_2$, observe that $\vert {y_2}_\alpha\vert\leq 1$ for every $\alpha\in\N^n_{2d}$, and 
from $M_d(\y_2-\y_1)\succeq0$, the diagonal elements ${y_2}_{2\alpha}-{y_1}_{2\alpha}$ are nonnegative.
Hence ${y_1}_{2\alpha}\leq {y_2}_{2\alpha}$ for every $\alpha\in\N^n_d$ and therefore,
\[\max\,[\,{y_1}_0,\,\max_{i=1,\ldots,n}L_{\y_1}(x_i^{2d})\,]\,\leq 1.\]
By \cite[Lemma 1]{lasserre-archiv},
this in turn implies that $\vert {y_1}_\alpha\vert\leq 1$ for every $\alpha\in\N^n_{2d}$, and so
the feasible set of $\Q_d$ is closed, bounded, hence compact, which in turn implies that
$\Q_d$ is solvable (i.e., has an optimal solution).

Let ${\y_1}^d$ be an optimal solution of $\Q_d$ and by completing with zeros,
make ${\y_1}^d$ an element of the unit ball $\B_\infty$ of $l_\infty$
(the Banach space of bounded sequences, equipped with the sup-norm). As $l_\infty$ is the topological dual of $l_1$, by the Banach-Alaoglu Theorem, $\B_\infty$ is
weak $\star$ compact, and even weak $\star$ sequentially compact; 
see e.g. Ash \cite{ash}. Therefore, there exists
 ${\y_1}^*\in \B_\infty$ and a subsequence $\{d_k\}\subset \N$ such that
 ${\y_1}^{d_k}\to {\y_1}^*$ as $k\to\infty$, for the weak $\star$ topology $\sigma(l_\infty,l_1)$. In particular, 
 \begin{equation}
 \label{pointwise}
\lim_{k\to\infty} {y_1}^{d_k}_\alpha\:=\: {y_1}^*_\alpha,\qquad \forall\,\alpha\in\N^n.
\end{equation}
Next let $d\in\N$ be fixed, arbitrary. From the pointwise convergence (\ref{pointwise})
we also obtain $M_{d}({\y_1}^*)\succeq0$ and $M_d(\y_2-{\y_1}^*)\succeq0$.
Similarly, $M_{d-r_j}(g_j{\y_1}^*)\succeq0$ for every $j=1,\ldots,m$.
As $d$ was arbitrary, by Theorem \ref{thput}, ${\y_1}^*$ has a representing
measure $\mu_1$ supported on $\K\subset\B$.
In particular, from (\ref{pointwise}), as $k\to\infty$,
\[\max\P\,\leq\,\max\Q_{d_k}\,=\,L_{\y_1^{d_k}}(p)\downarrow L_{\y_1^*}(p)\,=\,
\int pd\mu_1.\]
Next, as both $\mu_1$ and $\mu_2$ are supported on $[-1,1]^n$, and
$M_d(\y_2-{\y_1}^*)\succeq0$ for every $d$, one has
$\vert y_{2\alpha}-y^*_{1\alpha}\vert\leq 1$ for every $\alpha\in\N^n$. Hence
$\y_2-\y_1^*$ has a representing measure on $[-1,1]^n$.
As in the proof of Lemma \ref{ll}\footnote{If $\K\subset [-1,1]^n$ then in Lemma \ref{ll}, 
the condition $M_d(\y_2-\y_1)\succeq0,\:\forall d\in\N$, is sufficient.}, we conclude that
$\mu_1\leq\mu_2$. Therefore $\mu_1$ is admissible for problem $\P$, with
value $L_{\y_1^*}(p)=\int pd\mu_1\geq\max\P$. Therefore, $\mu_1$ must be an optimal solution of $\P$
(hence unique) and by Theorem \ref{th-mesures}, $L_{\y_1^*}=\int pd\mu_1=\int_\K pd\mu_2$. 
As the converging subsequence $\{d_k\}$ was arbitrary, it follows
that in fact the whole sequence ${\y_1}^d$ converges to ${\y_1}^*$ for the weak $\star$ topology
$\sigma(l_\infty,l_1)$. And so (\ref{conv}) holds.
This proves (a) and (b).
\end{proof}

Writing $M_d(\y_1)=\sum_\alpha A_\alpha {y_1}_\alpha$, and
$M_{d-r_j}(g_j\,\y_1)=\sum_\alpha B^j_\alpha {y_1}_\alpha$
for appropriate real symmetric matrices $(A_\alpha,B^j_\alpha)$,
the dual of $\Q_d$ reads:
\[\Q^*_d:\quad\left\{\begin{array}{ll}
\displaystyle\inf_{X,Y,Z_j} &\langle M_d(\y_2),Y\rangle\\
\mbox{s.t.}&\langle A_\alpha,Y-X\rangle-\displaystyle\sum_{j=1}^m\langle B^j_\alpha,Z_j\rangle=
p_{\alpha}\\
&X,Y,Z_j\succeq\,0,
\end{array}\right.\]
where $\langle X,Y\rangle=\mathrm{trace}\,(XY)$ is the standard inner product of real symmetric matrices, and $X\succeq0$ stands for $X$ is positive semidefinite. This
can be reformulated as:
\begin{equation}
\label{relaxdefdual}
\Q^*_d:\quad\left\{\begin{array}{ll}
\displaystyle\inf_{h,\sigma_0,\ldots,\sigma_m} &\displaystyle\int h\,d\mu_2\\
\mbox{s.t.}&h-p=\sigma_0+\displaystyle\sum_{j=1}^m\sigma_j\,g_j\\
&h\in\Sigma^2[x]_{d},\:\sigma_0\in\Sigma^2[x]_{d},\:\sigma_j\in\Sigma^2[x]_{d-r_j}.
\end{array}\right.
\end{equation}
The constraint of this semidefinite program
states that the polynomial $h-p$
is written in Putinar's form (\ref{putinarrep}) and so $h-p\geq0$ on $\K$.
In addition, $h\geq0$ because it is a sum of squares.

This interpretation of $\Q^*_d$ also shows why computing $\vol(\K)$ is difficult.
Indeed, when $p\equiv 1$, to get a good upper bound on $\vol(\K)$, one needs
to obtain a good {\it polynomial} approximation $h\in\R[x]$ of the indicator function
$I_\K(x)$ on $\B$. In general, high degree of $h$
will be necessary to attenuate side effects on the boundary of $\B$ and $\K$,
a well-known issue in interpolation with polynomials.
\begin{prop}
If $\K$ and $\B\setminus\K$ have a nonempty interior,
there is no duality gap, that is,
both optimal values of $\Q_d$ and $\Q^*_d$ are equal. In addition,
$\Q^*_d$ has an optimal solution 
$(h^*,(\sigma_j^*))$.
\end{prop}
\begin{proof}
Let $\mu_1$ be the uniform distribution
on $\K$, i.e., the restriction of $\mu_2$ to $\K$, and
let ${\y_1}=({y_1}_\alpha)$ be its sequence of moments up to degree $2d$.
As $\K$ has nonempty interior, then 
clearly $M_d(\y_1)\succ0$ and $M_{d-r_j}(g_j\,\y_1)\succ0$ for every
$j=1,\ldots,m$. If $\B\setminus\K$ also has nonempty interior then
$M_d(\y_2-\y_1)\succ0$ because
with $f\in\R[x]_d$ with coefficient vector $\f$,
\[\langle \f,M_d(\y_2-\y_1)\f\rangle\,=\,\int_{\B\setminus\K} f(x)^2d\mu_2,\quad\forall f\in\R[x]_d.\]
Therefore Slater's condition holds for
$\Q_d$ and the result follows from a standard result of duality in semidefinite programming; see e.g. \cite{boyd}.
\end{proof}
\begin{rem}
\label{rem1}
Let $f\in\R[x]$ and suppose that one wants to approximate the integral
$J^*:=\int_\K fd\mu_2$. Then for $d$ sufficiently large, an optimal solution of
$\Q_d$ allows to approximate $J^*$. Indeed,
\[J^*=\int_\K fd\mu_2\,=\,\int f\,d\mu_1\,=\,L_{{\y_1}^*}(f)\,=\,\sum_{\alpha\in\N^n}f_\alpha {y_1}^*_\alpha,\]
where ${\y_1}^*$ is the moment sequence of $\mu_1$, the
unique optimal solution of $\P$ (the restriction of $\mu_2$ to $\K$).
And so, from (\ref{conv}),
$L_{\y_1^d}(f)\approx J^*$ when $d$ is sufficiently large.
\end{rem}

\subsection{A hierarchy of linear programs}
\label{lp}
Let $\K\subset\B\subseteq [-1,1]^n$ be as in (\ref{setk}) and assume for simplicity that the $g_j$s have been scaled to satisfy
$0\leq g_j\leq 1$ on $\K$ for every $j=1,\ldots,m$.  In addition, assume that the family of polynomials
$(1,g_1,\ldots,g_m)$ generates the algebra $\R[x]$.   
For $d\in\N$, consider the following linear program:
\begin{equation}
\label{lprelax}
\L_d:\quad \left\{\begin{array}{ll}
\sup_{\y_1}&{y_1}_{0}\\
\mbox{s.t.}&L_{\y_2-\y_1}\left(\displaystyle\prod_{i=1}^n(1+x_i)^{\alpha_i}(1-x_i)^{\beta_i}\right)\,\geq\,0,\qquad
\alpha,\beta\in\N^n_d\\
&\\
&L_{\y_1}(g^\alpha\,(1-g)^\beta)\,\geq\,0,\qquad \alpha,\beta\in\N^n_d\end{array}\right.
\end{equation}
with optimal value denoted $\sup\L_d$ (and $\max\L_d$ if $\sup\L_d$ is finite). Notice that
$\sup\L_d\geq\vol(\K)/2^n$ for all $d$.
Indeed, the sequence $\y_1^*$ of moments of 
the Borel measure $\mu_1^*$ (restriction of $\mu_2$ to $\K$ and unique optimal solution of $\P$)
is a feasible solution of $\L_d$ for every $d$.

\begin{thm}\label{th-lp}
For the hierarchy of
linear programs $(\L_d)$ in (\ref{lprelax}), the following holds:\\

{\rm (a)} $\L_d$ has an optimal solution (i.e. $\sup\L_d=\max\L_d$)
and $\max\L_d\downarrow \vol(\K)/2^n$ as $d\to\infty$.

{\rm (b)} Let ${\y_1}^d$ be an optimal solution of $\L_d$. Then (\ref{conv}) holds.
\end{thm}
\begin{proof}
We first prove that $\L_d$ has finite value.
$\L_d$ always has a feasible solution $\y_1$, namely the moment vector associated with the Borel measure
$\mu_1$, the restriction of $\mu_2$ to $\K$, and so $\sup\L_d\geq\vol(\K)/2^n$. Next, from the constraint
$L_{\y_2-\y_1}(\bullet)\geq0$ with $\alpha=\beta=0$, we obtain ${y_1}_0\leq {y_2}_0\leq 1$.
Hence $\sup\L_d\leq 1$ and therefore, the linear program
$\L_d$ has an optimal solution ${\y_1}^d$. Fix $\gamma\in\N^n$ and $\epsilon>0$, arbitrary.
As $\vert x^\gamma\vert \leq 1<1+\epsilon$ on $\B$ (hence on $\K$), by Theorem \ref{thkrivine}(a),
\[1+\epsilon\pm x^\gamma=\sum_{\alpha,\beta\in\N^m}\,c^\gamma_{\alpha\beta}\,g^\alpha\,(1-g)^\beta,\]
for some $(c^\gamma_{\alpha\beta})\subset\R_+$ with $\vert\alpha\vert,\vert\beta\vert\leq s_\gamma$. 
Hence, as soon as $d\geq s_\gamma$, applying $L_{{\y_1}^d}$ yields
\[(1+\epsilon)\,{y_1}^d_0\pm {y_1}^d_\gamma\,=\,
\sum_{\alpha,\beta\in\N^m}\,c^\gamma_{\alpha\beta}\,
L_{\y_1}\left(g^\alpha\,(1-g)^\beta\right)\,\geq\,0,\]
and so 
\begin{equation}
\label{bound}
\forall \gamma\in\N^n:\qquad \vert {y_1}^d_\gamma\vert\leq (1+\epsilon)\,{y_1}_0^d\leq 1+\epsilon,\qquad \forall d\geq s_\gamma.
\end{equation}
Complete ${\y_1}^d$ with zeros to make it an element of $\R^\infty$. Because of (\ref{bound}), using a standard diagonal element, there exists a subsequence $(d_k)$ and an element ${\y_1}^*\in (1+\epsilon)\,\B_\infty$ (where $\B_\infty$ is the unit ball of $l_\infty$) 
such that (\ref{pointwise}) holds.
Now with $\alpha,\beta\in \N^m$ fixed, arbitrary, (\ref{pointwise}) yields
$L_{{\y_1}^*}(g^\alpha\,(1-g)^\beta)\geq0$. Hence by Theorem \ref{thkrivine}(b),
${\y_1}^*$ has a representing measure $\mu_1$ supported on $\K$.
Next, let $\y_0:=\y_2-{\y_1}^*$. Again, (\ref{pointwise}) yields:
\[L_{\y_0}\left(\prod_{i=1}^n(1+x_i)^{\alpha_i}\,(1-x_i)^{\beta_i}\right)\,\geq\,0,\qquad \forall\alpha,\beta\in\N^n,\]
and so by Theorem \ref{thkrivine}(b), $\y_0$ is the moment vector of some Borel measure 
$\mu_0$ supported on $[-1,1]^n$. As measures on compact sets are identified with their moments,
and ${y_0}_\alpha+{y_1}^*_\alpha={y_2}_\alpha$ for every $\alpha\in\N^n$,
it follows that $\mu_0+\mu_1=\mu_2$, and so $\mu_1\leq\mu_2$. Therefore, $\mu_1$ is an admissible solution to $\P$ with parameter $p\equiv 1$, and with value $\mu_1(\K)={y_1}^*_0\geq\vol(\K)/2^n$. Hence, $\mu_1$ is the unique optimal solution to $\P$ with value 
$\mu_1(\K)=\vol(\K)/2^n$.

Finally, by using (\ref{pointwise}) and following the same argument
as in the proof of Theorem \ref{th-sdp}, one obtains the desired result
 (\ref{conv}).
\end{proof}

Remark \ref{rem1} also applies to the LP relaxations (\ref{lprelax}).

\subsection{Integration against a weight function}

With $\K\subset\B$ as in (\ref{setk}) suppose now that one wishes to approximate the integral
\begin{equation}
\label{weight}
J^*\,:=\,\int_\K f(x)\,w(x)\,dx,\end{equation}
for some given nonnegative {\it weight} function $w:\R^n\to \R$,
and where $f\in\R[x]_d$ is some nonnegative
polynomial.  One makes the following assumption:

\begin{assumption}
\label{ass2}
One knows the moments $\y_2=({y_2}_\alpha)$ of the Borel measure $d\mu_2=wdx$ on $\B$, that is:
\begin{equation}
\label{mom}
{y_2}_\alpha\,=\,\int_{\B}x^\alpha \,d\mu_2\:\left(=\,\int_{\B}x^\alpha\, w(x)\,dx\,\right),\qquad \alpha\in\N^n.
\end{equation}
\end{assumption}
Indeed, for many weight functions $w$, and given $d\in\N$, one may compute
the moments $\y_2=({y_2}_\alpha)$ of $\mu_2$
via {\it cubature} formula, exact up to degree $d$. In practice, one only knows 
finitely many moments of $\mu_2$, say up to degree $d$, fixed.

Consider the hierarchy of semidefinite programs 
\begin{equation}
\label{relax-inte}
\Q_d:\quad\left\{\begin{array}{lll}
\displaystyle\sup_{\y_1} &L_{\y_1}(f)&\\
\mbox{s.t.}&M_d(\y_1)&\succeq\,0\\
&M_d(\y_2-\y_1)&\succeq\,0\\
&M_{d-r_j}(g_j\,\y_1)&\succeq\,0,\qquad j=1,\ldots,m
\end{array}\right.
\end{equation}
with $\y_2$ as in Assumption \ref{ass2}.
\begin{thm}
\label{inte}
Let Assumption \ref{assput} and \ref{ass2}  hold and consider the hierarchy of
semidefinite programs $(\Q_d)$ in (\ref{relax-inte}) with $\y_2$ as in (\ref{mom}). 
Then $\Q_d$ is solvable and $\max\Q_d\downarrow J^*$ as $d\to\infty$.
\end{thm}

The proof is almost a verbatim copy of that of Theorem \ref{th-sdp}.

\section{Numerical experiments and discussion}

In this section we report some numerical experiments
carried out with Matlab and the package GloptiPoly 3 for manipulating and
solving generalized problems of moments \cite{gloptipoly3}.
The SDP problems were solved with SeDuMi 1.1R3 \cite{sedumi}.
Univariate Chebyshev polynomials were manipulated with the
{\tt chebfun} package \cite{chebfun}.

The single-interval example below permit to visualize the 
numerical behavior of the algorithm.
The folium example
illustrates that, as expected,  the non-convexity of $\K$ does not
seem to penalize the moment approach.
Finally, our experience reveals that the choice of 
alternative polynomial bases affects the quality
of the approximations.

\subsection{Single interval}

Consider the elementary one-dimensional set $\K = [0,\:\frac{1}{2}] = 
\{x \in {\mathbb R} \: :\: g_1(x) = x(\frac{1}{2}-x) \geq 0\}$ included in the
unit interval $\B=[-1,1]$. We want to
approximate $\vol(\K) = \frac{1}{2}$.
Moments of the Lebesgue measure $\mu_2$ on $\B$ are given by
$\y_2 = (2,0,2/3,0,2/5,0,2/7,\ldots)$.

Here is a simple Matlab script using GloptiPoly 3 instructions
to input and solve the SDP relaxation $\Q_d$ of the LP moment problem $\P$ with $p\equiv 1$:
\begin{verbatim}
>> d = 10; % degree
>> mpol x0 x1
>> m0 = meas(x0); m1 = meas(x1);
>> g1 = x1*(1/2-x1);
>> dm = (1+(0:d))'; y2 = ((+1).^dm-(-1).^dm)./dm; 
>> y0 = mom(mmon(x0,d)); y1 = mom(mmon(x1,d)); 
>> P = msdp(max(mass(m1)), g1>=0, y0==y2-y1); % input moment problem
>> msol(P); % solve SDP relaxation
>> y1 = double(mvec(m1)); % retrieve moment vector
\end{verbatim}
The volume estimate is then the first entry in vector $\tt y1$.
Note in particular the use of the moment constraint {\tt y0==y2-y1} which
ensures that moments $\y_0$ of $\mu_0$ will be substituted
by linear combinations of moments $\y_1$ of $\mu_1$ (decision
variables) and moments $\y_2$ of $\mu_2$ (given).

\begin{figure}[h!]
\begin{center}
\includegraphics[width=\figurewidth\linewidth]{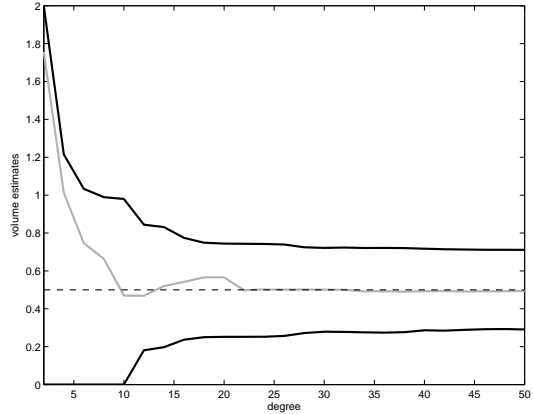}
\caption{Three sequences of approximations of $\vol[0,\frac{1}{2}]$ obtained by
solving SDP relaxations of increasing degree.\label{oneinterval}}
\end{center}
\end{figure}

Figure \ref{oneinterval} displays three approximation sequences
of $\vol(\K)$ obtained by solving SDP
relaxations (\ref{relax-inte}) of increasing degrees $d=2,\ldots,50$
of the infinite-dimensional LP moment problem $\P$ with three different parameters $p$:
\begin{itemize}
\item the upper curve (in black) is a monotone non increasing sequence
of upper bounds obtained by maximizing $\int d\mu_1$, the mass of $\mu_1$,
using the objective function {\tt max(mass(m1))} in the above script;
\item the medium curve (in gray) is a sequence of approximations
obtained by maximizing $\int pd\mu_1$ with $p:=g_1$, using the objective function
$\tt max(g1)$ in the above script;
\item the lower curve (in black) is a monotone non decreasing sequence
of lower bounds on $\vol(\K)$ obtained by computing 
upper bounds on the volume of $\B\setminus\K$, using
the objective function {\tt max(mass(m1))} and the support
constraint {\tt g1<=0} in the above script. The volume estimate
is then {\tt 2-y1(1)}.
\end{itemize}
We observe a much faster convergence when maximizing $\int g_1d\mu_1$
instead of $\int d\mu_1$; the upper and
lower curves apparently exhibit slow convergence.

\begin{figure}[h!]
\begin{center}
\includegraphics[width=\figurewidth\linewidth]{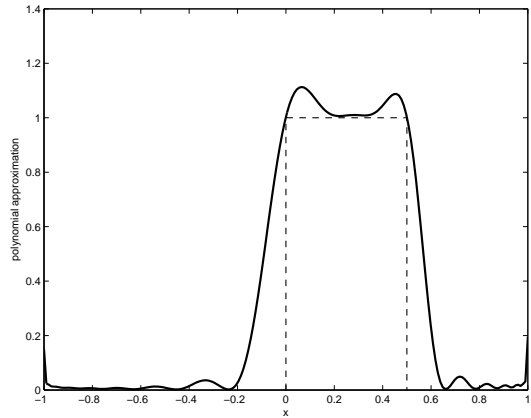}
\caption{Positive polynomial approximation of degree 50 (solid)
of the indicator function $I_{[0,\frac{1}{2}]}$ (dashed) on $[-1,1]$.
\label{oneinterval_up}}
\end{center}
\end{figure}

\begin{figure}[h!]
\begin{center}
\includegraphics[width=\figurewidth\linewidth]{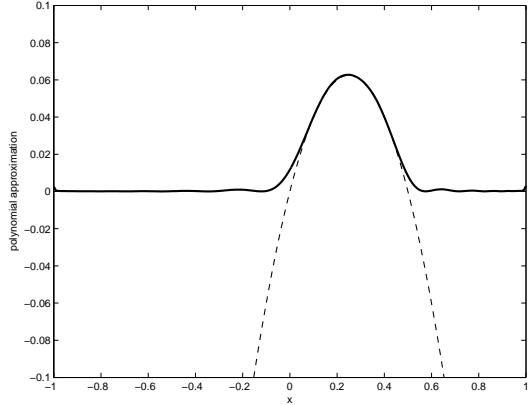}
\caption{Positive polynomial approximation of degree 50 (solid) of the
positive piecewise-polynomial function $\max(0,g_1)$ on $[-1,1]$.
Polynomial $g_1$ is represented in dashed line.
\label{oneinterval_poly}}
\end{center}
\end{figure}

\begin{figure}[h!]
\begin{center}
\includegraphics[width=\figurewidth\linewidth]{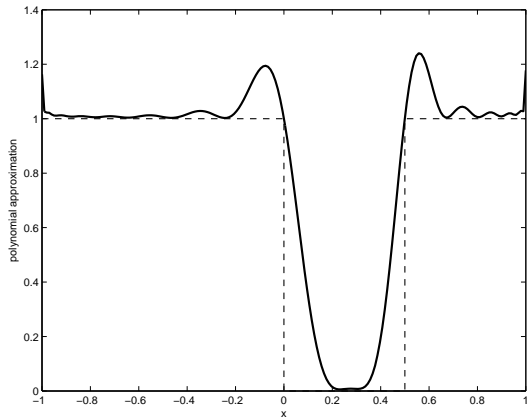}
\caption{Positive polynomial approximation of degree 50 (solid) of the
complementary indicator function $1-I_{[0,\frac{1}{2}]}$ (dashed) on $[-1,1]$.
\label{oneinterval_down}}
\end{center}
\end{figure}

To analyze these phenomena, we use solutions of the dual
SDP problems, provided automatically by the primal-dual
interior-point method implemented in the SDP solver SeDuMi.
On Figure \ref{oneinterval_up} we represent the degree-50 positive polynomial
approximation $h$ of the
indicator function $I_{\K}$ on $\B$, which minimizes $\int_\B h dx$ while
satisfying $h-1 \geq 0$ on $\K$ and $h \geq 0$ on $\B\setminus\K$ (yielding
the volume estimate of the upper curve in Figure \ref{oneinterval}).
On Figure
\ref{oneinterval_poly}, we represent the degree-50 polynomial approximation
$h$ of the piecewise-polynomial function
$\max(0,g_1)$, which minimizes $\int_\B h dx$ while satisfying
$h-g_1 \geq 0$ on $\K$ and $h \geq 0$ on $\B\setminus\K$ (yielding the
volume estimate of 
the medium curve in Figure \ref{oneinterval}).
On Figure \ref{oneinterval_down}
we represent the degree-50 polynomial approximation $h$ of
the complementary indicator function $1-I_{\K}$, which minimizes
$\int_\B h dx$ while satisfying $h-1 \geq 0$ on $\B\setminus\K$ and
$h \geq 0$ on $\K$
(yielding the volume estimate of
the lower curve in Figure \ref{oneinterval}).
We observe the characteristic oscillation phenomena
near the boundary, typical of polynomial approximation
problems \cite{trefethen}. The continuous function $\max(0,g_1)$
is easier to approximate than discontinuous indicator functions, and this partly explains
the better convergence of the medium approximation
on Figure \ref{oneinterval}.

On Figures \ref{oneinterval_up} and \ref{oneinterval_down},
one observes relatively large oscillations near the boundary points
$x\in\{-1,0,\frac{1}{2},1\}$ which significantly corrupt
the quality of the volume approximation. To some extent, these oscillations
can be reduced by using a Chebyshev polynomial
basis instead of the standard power basis. 

\begin{figure}[h!]
\begin{center}
\includegraphics[width=\figurewidth\linewidth]{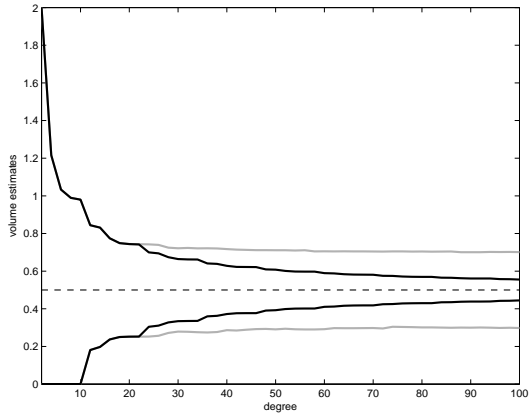}
\caption{Upper and lower bounds on $\vol[0,\frac{1}{2}]$ obtained by
solving SDP relaxations in the Chebyshev basis (black) and power basis
(gray).\label{oneintervalcheb}}
\end{center}
\end{figure}

Figure \ref{oneintervalcheb} displays upper and lower
bounds on the volume, computed up to degree 100, with
the power basis (in gray) and with the Chebyshev basis (in black).
Note that in order to input and solve SDP problems in the Chebyshev
basis, we used our own implementation and the {\tt chebfun} package
since GloptiPoly 3 supports only the power basis. In Figure
\ref{oneintervalcheb} we see that above degree 20 the quality
of the bounds obtained with the power basis deteriorates, which suggests
that the SDP solver encounters some numerical problems
rather than convergence becoming slower (which is confirmed 
when changing to Chebyshev basis; see below).
It seems that the SDP solver is not able to improve the
bounds, most likely due to the symmetric Hankel structure
of the moment matrices in the power basis: indeed, it is known
that the conditioning (ratio of extreme singular values) of positive 
definite Hankel matrices is an exponential function of
the matrix size \cite{higham}. When the smallest singular values
reach machine precision, the SDP solver is not able to optimize
the objective function any further.

\begin{figure}[h!]
\begin{center}
\includegraphics[width=\figurewidth\linewidth]{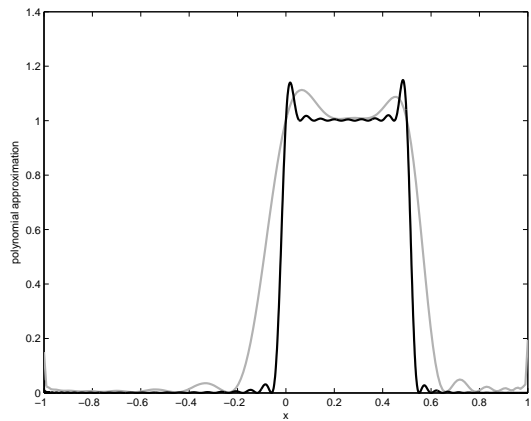}
\caption{Positive polynomial approximation of degree 100 
of the indicator function $I_{[0,\frac{1}{2}]}$ in the
Chebyshev basis (black) and power basis (gray).
\label{oneintervalcheb_up}}
\end{center}
\end{figure}

\begin{figure}[h!]
\begin{center}
\includegraphics[width=\figurewidth\linewidth]{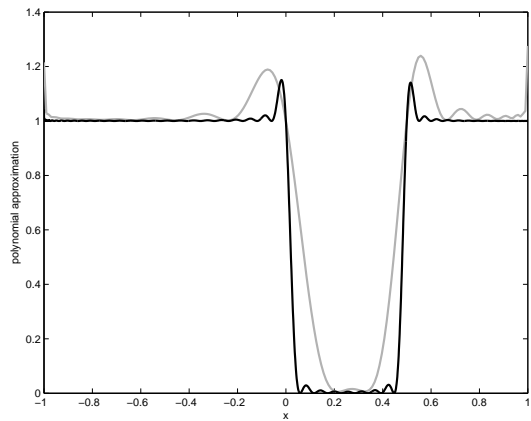}
\caption{Positive polynomial approximation of degree 100 
of the indicator function $1-I_{[0,\frac{1}{2}]}$ in the
Chebyshev basis (black) and power basis (gray).
\label{oneintervalcheb_down}}
\end{center}
\end{figure}

In Figures \ref{oneintervalcheb_up} and \ref{oneintervalcheb_down}
one observes that the degree-100 polynomial approximation $h(x)$
of the indicator function and its complement
are tighter in the Chebyshev basis (black) than in the
power basis (gray). Firstly, we observe that the degree-100
approximations in the power basis do not significantly differ
from the degree-50 approximations in the same basis,
represented in Figures \ref{oneinterval_up} and \ref{oneinterval_down}.
This is consistent with the very flat behavior of the right half of
the upper and lower curves (in gray) in Figure \ref{oneintervalcheb}.
Secondly, some coefficients of $h(x)$ in the power basis
have large magnitude $h(x) = 1.0019+3.6161x-29.948x^2+
\cdots +88123x^{49}+54985x^{50}+\cdots-1018.4x^{99}+26669x^{100}$
with the Euclidean norm of the coefficient vector greater
than $10^6$. In contrast, the polynomial $h(x)$ obtained
in the Chebyshev basis $h(x) = 0.1862t_0(x)+0.093432t_1(x)
-0.30222t_2(x)+\cdots +0.0055367t_{49}(x)-0.020488t_{50}(x)+\cdots
-0.0012267t_{99}(x)+0.0011190t_{100}(x)$
has a coefficient vector of Euclidean norm around $0.57627$,
where $t_k(x)$ denotes the $k$-th Chebyshev polynomial.
Thirdly, oscillations around points $x=0$ and $x=1/2$
did not disappear with the Chebyshev basis,
but the peaks are much thinner than with the power basis.
Finally, the oscillations near the interval ends $x=-1$
and $x=1$ are almost suppressed, a well-known property of
Chebyshev polynomials which have a denser root distribution
near the interval ends. 

From these simple observations, we 
conjecture that a polynomial basis with a dense root
distribution near the boundary of the semi-algebraic sets $\K$
and $\B$ should ensure a better convergence of the hierarchy
of volume estimates.

\begin{figure}[h!]
\begin{center}
\includegraphics[width=\figurewidth\linewidth]{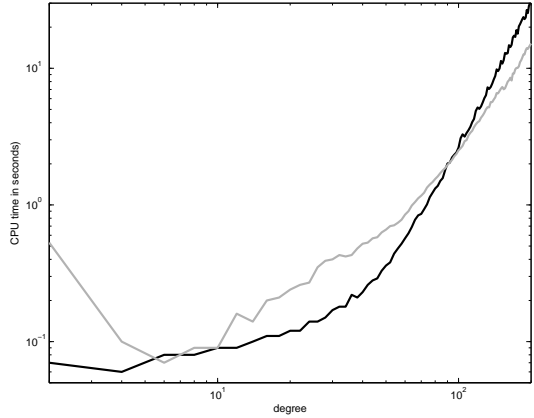}
\caption{CPU time required to solve the SDP relaxations (Chebyshev
basis in black, power basis in gray) as a function of the degree. 
\label{timecheb}}
\end{center}
\end{figure}

Finally, Figure \ref{timecheb} displays the CPU time
required to solve the SDP problems (with SeDuMi, in the power basis in gray and
in the Chebyshev basis in black) as a function of the degree,
showing a polynomial dependence slightly slower than cubic in the power basis
(due to the sparsity of moment matrices) and slightly faster than cubic
in the Chebyshev basis.
For example, solving the SDP problem of degree 100 takes
about 2.5 seconds of CPU time on our standard desktop computer.

%
%

\subsection{Bean}

\begin{figure}[h!]
\begin{center}
\includegraphics[width=\figurewidth\linewidth]{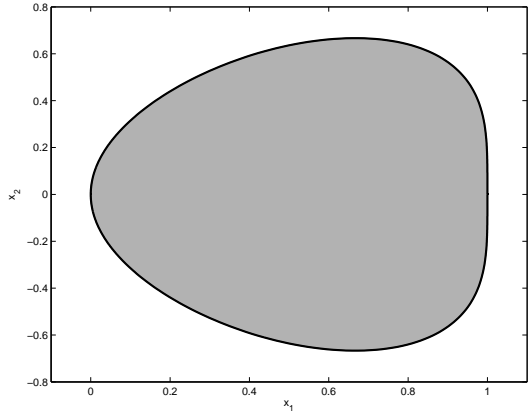}
\caption{Bean surface.\label{bean2d}}
\end{center}
\end{figure}

Consider $\K = \{x \in {\mathbb R}^2 \: :\: g_1(x) = x_1(x^2_1+x^2_2)-(x^4_1+x^2_1x^2_2+x^4_2) \geq 0
\}$ displayed in Figure \ref{bean2d}, which is a surface delimited by an algebraic curve $g_1(x) = 0$
of genus zero, hence rationally parametrizable. From the parametrization
$x_1(t) = (1+t^2)/(1+t^2+t^4)$, $x_2(t) = tx_1(t)$, $t \in \mathbb R$,
obtained with the {\tt algcurves}
package of {\tt Maple}, we can calculate
\[
\begin{array}{rcl}
\vol(\K)&=& \int_{\K} dx_1dx_2 = \int_\R x_1(t)dx_2(t) =
\int_{\R} \frac{(1-t)(1+t)(1+t^2)(1+3t^2+t^4)}{(1+t+t^2)^3(1-t+t^2)^3}dt \\
&=& \frac{7\sqrt{3}\pi}{36} \approx 1.0581
\end{array}
\]
with the help of the {\tt int} integration routine of {\tt Maple}.
Similarly, we can calculate symbolically the first moments of the Lebesgue measure $\mu_1$
on $\K$, namely ${y_1}_{00} = \vol(\K)$, ${y_1}_{10} = \frac{23}{42}\vol(\K)$,
${y_1}_{01} = 0$, ${y_1}_{20} = \frac{23}{63}\vol(\K)$, ${y_1}_{11} = 0$,
${y_1}_{02} = \frac{113}{1008}\vol(\K)$ etc.
Observe that $\K \subseteq \B = [-1,1]^2$.

\begin{figure}[h!]
\begin{center}
\includegraphics[width=\figurewidth\linewidth]{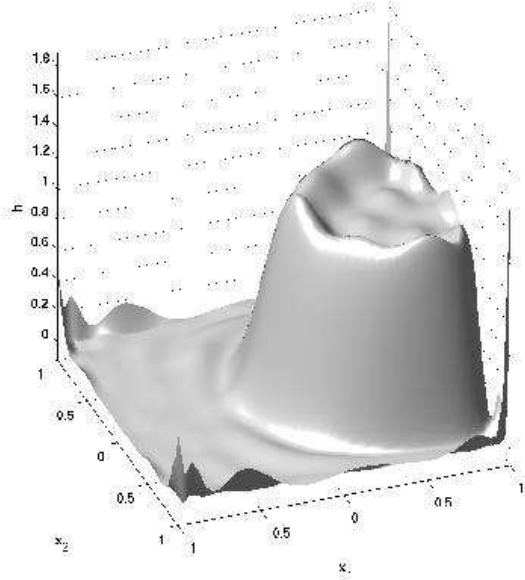}
\caption{Positive polynomial approximation of degree 20
of the indicator function of the bean surface.
\label{bean3d}}
\end{center}
\end{figure}

On Figure \ref{bean3d} we represent a degree-20 positive polynomial approximation $h$
of the indicator function $I_\K$ on $\B$
obtained by solving an SDP problem with 231 unknown moments.
We observe the typical oscillations near the boundary regions,
but we can recognize the shape of Figure \ref{bean2d}.

\begin{table}[h!]
\begin{center}
\begin{tabular}{c|ccccccccc}
degree & 2 & 4 & 6 & 8 & 10 & 12 & 14 \\ \hline
error & 78\% & 63\% & 13\% & 0.83\% & 9.1\% & 0.80\%  & 3.31\% \\[1em]
degree & 16 & 18 & 20 & 22 & 24 & 26 & 28 & 30 \\ \hline
error & 3.8\% & 3.3\% & 2.6\% & 5.6\%  & 4.1\%  & 4.1\% & 3.9\% & 3.7\% 
\end{tabular}\\[5mm]
\caption{Relative error when approximating the volume of the
bean surface, as a function of the degree of the SDP relaxation.\label{beanpoly}}
\end{center}
\end{table}

In Table \ref{beanpoly} we give relative errors in percentage
observed when solving successive SDP relaxations (in the power basis)
of the LP moment problems of maximizing $\int g_1d\mu_1$.
Note that the error sequence is not monotonically decreasing since we
do not maximize $\int d\mu_1$ and a good approximation can be obtained with few moments.
Above degree 16, the approximation stagnates around $4\%$,
Most likely this is due to the use of the power basis,
as already observed in the previous univariate examples.
For example, at degree 20, one obtains the 6 first
moment approximation
\[{y_2}_{00}^{20} = 1.10,\:{y_2}_{10}^{20}=0.589, \:{y_2}_{01}^{20}=0.00,
\:{y_2}_{20}^{20}=0.390,\:{y_2}_{11}^{20}=0.00,\:{y_2}_{02}^{20}=0.122\]
to be compared with the exact numerical
values 
\[{y_2}_{00} = 1.06, \:{y_2}_{10}=0.579,\:{y_2}_{01}=0.00,
\:{y_2}_{20}=0.386,\:{y_2}_{11}=0.00,\:{y_2}_{02}=0.119.\]
Increasing the degree does not provide a better approximation.
It is expected that a change of basis (e.g. multivariate Chebyshev
or trigonometric)
can be useful in this context.

\subsection{Folium}

\begin{figure}[h!]
\begin{center}
\includegraphics[width=\figurewidth\linewidth]{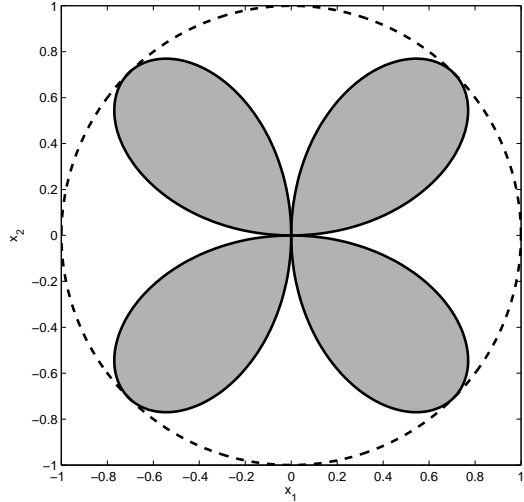}
\caption{Folium surface. 
\label{folium2d}}
\end{center}
\end{figure}

Consider $\K = \{x \in {\mathbb R}^2 \: :\: g_1(x) = -(x^2_1+x^2_2)^3+4x^2_1x^2_2 \geq 0
\}$ displayed in Figure \ref{folium2d},
which is a surface delimited by an algebraic curve of polar equation
$\rho = \sin(2\theta)$. The surface is contained in the unit disk $\B$,
on which the Lebesgue measure has moments
\[
{y_2}_{\alpha} = \frac{(1+(-1)^{\alpha_1})(1+(-1)^{\alpha_2})
\Gamma(\frac{1}{2}(1+\alpha_1))\Gamma(\frac{1}{2}(1+\alpha_2))}
{\Gamma(\frac{1}{2}(4+\alpha_1+\alpha_2))},\qquad\forall\,\alpha\in\N^2,
\]
where $\Gamma$ denotes the gamma function.
The area is $\vol(\K) = \frac{1}{2}\int_0^{2\pi}\sin^2(2\theta)d\theta =
\frac{1}{2}\pi$ and so, $\vol(\K\setminus\B) = \pi-\vol(\K) = \frac{1}{2}\pi$.

In Table \ref{folium} we give relative errors in percentage
observed when solving successive SDP relaxations (in the power basis)
of the LP moment problems of maximizing $\int g_1d\mu_1$.
We observe that nonconvexity of $\K$ does not play any special
role. The quality of estimates does not really improve for degrees greater than 20.
Here too, an alternative polynomial basis with dense root distribution
near the boundaries of $\K$ and $\B$ would certainly help.

\begin{table}[h!]
\begin{center}
\begin{tabular}{c|ccccccccc}
degree & 4 & 6 & 8 & 10 & 12 & 14 & 16 \\ \hline
error & 87\% & 19\% & 14\% & 9.4\% & 4.3\% & 4.5\%  & 5.9\% \\[1em]
degree & 18 & 20 & 22 & 24 & 26 & 28 & 30 \\ \hline
error & 1.2\% & 5.3\% & 5.9\% & 7.2\%  & 8.7\%  & 9.0\% & 8.8\%
\end{tabular}\\[5mm]
\caption{Relative error when approximating the volume of the
folium surface, as a function of the degree of the SDP relaxation.\label{folium}}
\end{center}
\end{table}



\begin{figure}[h!]
\begin{center}
\includegraphics[width=\figurewidth\linewidth]{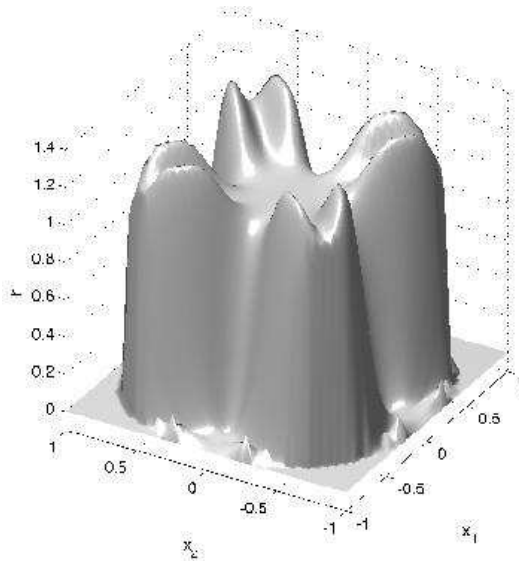}
\caption{Positive polynomial approximation of degree 20
of the indicator function of the folium surface.
\label{folium3d}}
\end{center}
\end{figure}

Figure \ref{folium3d} displays a degree-20 positive polynomial approximation $h$
of the indicator function $I_\K$ on $\B$
obtained by solving an SDP problem with 231 unknown moments.
For visualization purposes, $\max(5/4,h)$ rather than $h$ is displayed.
Again typical oscillations occur near the boundary regions,
but we can recognize the shape of Figure \ref{folium2d}.

\section{Concluding Remarks}

The methodology presented in this paper is general enough and applies to
compact basic semi-algebraic sets which are neither necessarily convex nor connected.
Its efficiency is related to the 
degree needed to obtain a good polynomial approximation
of the indicator function of $\K$ (on a simple set that contains $\K$) and 
from this viewpoint,  convexity of $\K$ does not help much. 
On the other hand, the method is limited by the size of problems that
SDP solvers presently available can handle.
Moreover, the impact of the choice of the polynomial basis (e.g., Chebyshev or
trigonometric) on the quality of the solution of the SDP relaxations
deserves further investigation for a better understanding.
Therefore, in view of the present status of SDP solvers 
and since in general high accuracy will require high degree,
the method can provide good approximations for
problems of small dimension (typically $n=2$ or $n=3$).
However, if one is satisfied with cruder bounds then one may consider problems in higher dimensions.

\end{document}